\documentclass[a4paper,11pt]{scrartcl}
\usepackage{amsmath,amssymb}
\usepackage{amsthm}
\usepackage{graphicx}
\usepackage[square,sort&compress,numbers]{natbib}
\usepackage{url}
\usepackage{type1cm}
\usepackage{mathtools}\mathtoolsset{showonlyrefs=false}

\usepackage{subcaption}
\newcommand{\Authornote}{\renewcommand{\thefootnote}{\fnsymbol{footnote}}}
\newcommand{\authornote}{\Authornote\footnote}
\theoremstyle{plain}
\newtheorem{theorem}{Theorem}[section]

\theoremstyle{definition}

\theoremstyle{remark}
\newtheorem{remark}[theorem]{Remark}

\newcommand{\refrem}[1]{Remark~\ref{#1}}
\newcommand{\reffig}[1]{Figure~\ref{#1}}

\newcommand{\finbox}{\nolinebreak\hfill{\small $\blacksquare$}}

\newcommand{\convex}{\mathop{\mathrm{co}}\nolimits}
\newcommand{\argmin}{\operatornamewithlimits{\mathrm{arg\,min}}}

\renewcommand{\Re}{\ensuremath{\mathbb{R}}}

\newcommand{\bi}[1]{\ensuremath{\boldsymbol{#1}}}

\newcommand{\NC}{\ensuremath{\mathcal{N}}}

\begin{document}


\begin{center}
  {\Large\bfseries\sffamily%
  Simple Heuristic for Data-Driven Computational Elasticity }\\
  \medskip
  {\Large\bfseries\sffamily%
  with Material Data Involving Noise and Outliers: }\\
  \medskip
  {\Large\bfseries\sffamily%
  A Local Robust Regression Approach 
  }%
  \par%
  \bigskip%
  {
  Yoshihiro Kanno~\authornote[2]{%
    Mathematics and Informatics Center, 
    The University of Tokyo, 
    Hongo 7-3-1, Tokyo 113-8656, Japan.
    E-mail: \texttt{kanno@mist.i.u-tokyo.ac.jp}. 
    }
  }
\end{center}

\begin{abstract}
  Data-driven computing in applied mechanics utilizes the material data set 
  directly, and hence is free from errors and uncertainties stemming from 
  the conventional material modeling. 
  This paper presents a simple heuristic for 
  data-driven computing, that is robust against 
  noise and outliers in a data set. 
  For each structural element, we extract the material property from 
  some nearest data points. Using the nearest neighbors reduces the 
  influence of noise, 
  compared with the existing method that uses a single data point. 
  Also, the robust regression is adopted to reduce the influence of 
  outliers. 
  Numerical experiments on the static equilibrium analysis of trusses are 
  performed to illustrate that the proposed method is robust against the 
  presence of noise and outliers and, hence, is effective for dealing 
  with real-world data. 
\end{abstract}

\begin{quote}
  \textbf{Keywords}
  \par
  Data-driven computing; 
  model-free computational mechanics; 
  outlier; 
  local regression; 
  robust statistics. 
\end{quote}

\section{Introduction}

Recent development and spread of data science are enormously 
remarkable \cite{Mat13,TLCV15,FH17,SZ16,CCS12,FATAKZFB14}. 
It has been widely recognized that the methodology of knowledge 
extraction from data is extremely useful in diverse fields. 

Database methods have also been developed in engineering computation. 
In the area of computer graphics, data-driven methods were proposed to 
construct simulation models of elastic fabrics \cite{WObR11,MBTBMOM12}, 
where cloth deformation models are estimated from data of 
experimental measurements. 
In computational mechanics, macroscopic constitutive properties of 
composites are extracted from a data set of the results of numerical 
material tests \cite{TKHIY13,KO15,TZ07,TW07,CSY12,BBLHABCL17}. 

This paper is inspired by the work of \citet{KO16}, 
where the paradigm of data-driven computational mechanics is presented. 
The awareness of issues raised in \cite{KO16} is summarized as follows. 
As an example, consider the static equilibrium analysis of a structure. 
The boundary value problem consists of 
(i) the compatibility relation, 
(ii) the force-balance equation, and 
(iii) the constitutive law. 
Among them, (i) is the kinematic constraint, 
and (ii) is derived from Newton's laws of motion. 
Therefore, (i) and (ii) do not possess any uncertainty or error. 
In contrast, (iii) is a formulation obtained through a physical modeling 
based on experiments, and hence is empirical and uncertain. 
Based on this observation, \citet{KO16} attempts to directly utilize a data set 
obtained from the physical experiments, without resorting to empirical 
modeling of the material in (iii). 
In the framework of this data-driving computing, we suppose that a 
material data set in the stress--strain space is given. 
Then we find a solution, that satisfies (i) and (ii) strictly 
(in the same manner as the conventional framework requires) and is 
closest to the data set in the stress--strain space. 
Thus, the data-driven computing does not require any modeling in (iii), 
and hence is free from errors and uncertainties stemming from the 
material modeling. 

The method proposed in 
\citep{KO16} has recently been extended to static problems 
with geometrical nonlinearity \citep{NK18}, 
three-dimensional continua \citep{CMO18}, 
and dynamic problems \citep{KO18}. 
Independently, another data-driven approach that makes use of the 
manifold learning for estimating the material law has been developed 
\cite{IAcAGCC18,IBAAcCLC17}. 

In \cite{KO16}, in the stress--strain space the distance between a point 
and a data set is defined as the distance between the point and the 
closest data point in the data set. 
This means that, for each structural element, the information of only 
one data point is adopted to extract the material property. 
Therefore, even if the data set consists of a large number of data points, 
most information that the data set has is not utilized for computation. 
Using only a single data point also means that the method has high 
estimation variance, i.e., it is sensitive to small fluctuation in the data set. 
In other words, the computational result will be seriously affected by 
noise and/or outliers in the data set. 

Attention of this paper is focused on the robustness of data-driven 
computational elasticity methods against noise and outliers involved 
in a data set. 
Specifically, instead of the single closest data point, 
we attempt to use the information of the $k$ nearest data 
points for each structural element, where $k$ is a sufficiently small 
positive constant compared with the total number of the data points in 
the data set. 
It is expected that using several data points reduces the influence of 
noise involved in the data set. 
To extract the material property, we apply a robust regression method to 
the $k$ nearest data points. 
This reduces the influence of outliers. 
Thus, the estimation variance of the data-driven solver is reduced. 
Also, locality of the regression avoids the bias of oversmoothing. 
Direct comparison of the computational results of the proposed method 
and the existing method in \cite{KO16} appears in section~\ref{sec:ex.1}. 

Recently, \citet{KO17} has proposed to incorporate a maximum-entropy 
estimation into their original data-driven solver \cite{KO16} in order 
to deal with a data set involving noise. 
In this method, the solution to be found is defined as a global optimal 
solution of a nonconvex optimization problem. 
As a heuristic, \citet{KO17} used simulated annealing to 
approximately solve this complicated nonconvex optimization problem. 
This paper attempts to present an alternative method that is much 
simpler. 
This method, based on the local robust regression, is also nothing 
more than a heuristic, in the sense that convergence is not guaranteed. 
However, it is simple, easy to understand, and easy to implement. 
This feature is considered attractive from the viewpoint of the quality 
management of numerical simulation.

The paper is organized as follows. 
In section~\ref{sec:overview}, we discuss influence of noise and outliers 
on a data-driven approach to computational elasticity, 
and give a brief overview of the method proposed in this paper. 
Section~\ref{sec:method} presents a full description of the proposed method. 
Section~\ref{sec:ex} reports the results of numerical experiments. 
Some conclusions are drawn in section~\ref{sec:conclude}.


In our notation, 
${}^{\top}$ denotes the transpose of a vector or a matrix. 
For simplicity, we often write the $(n+m)$-dimensional column vector 
$(\bi{x}^{\top} , \bi{y}^{\top})^{\top}$ consisting of 
$\bi{x} \in \Re^{n}$ and $\bi{y} \in \Re^{m}$ as 
$(\bi{x}, \bi{y})$. 
We use $\convex S$ to denote the convex hull of set $S$. 
For finite set $T$, the cardinality of $T$, i.e., 
the number of elements in $T$, is denoted by $|T|$. 
We use $\NC(\mu,\sigma^{2})$ to denote 
the normal distribution with mean $\mu$ and variance $\sigma^{2}$.

\section{Motivation and overview}
\label{sec:overview}

In this section, we provide motivation and brief overview of the 
method presented in this paper. 

\begin{figure}[tp]
  \centering
  \includegraphics[scale=0.70]{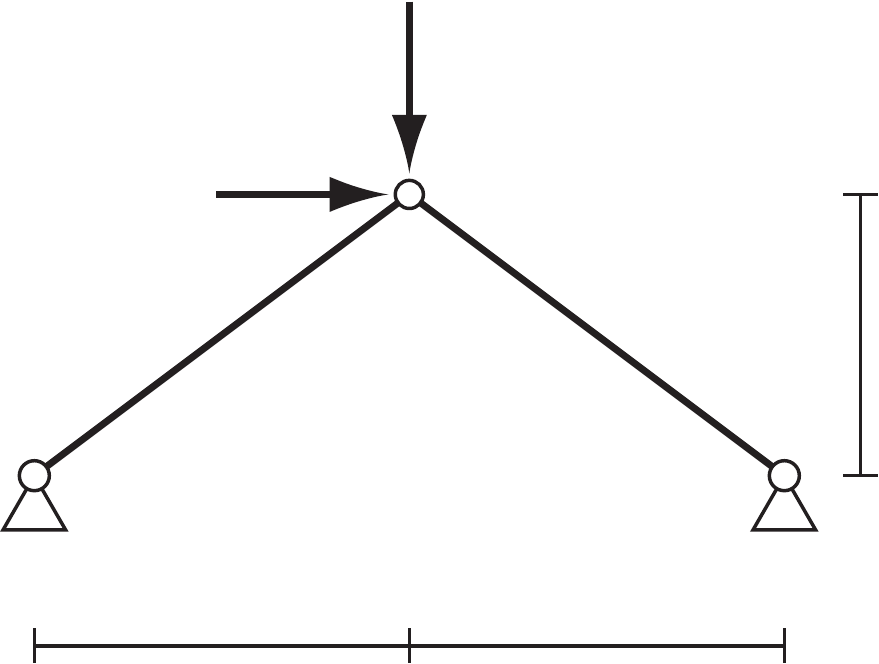}
  \begin{picture}(0,0)
    \put(-200,-50){
    \put(60,44){{\footnotesize $4\,\mathrm{m}$}}
    \put(133,44){{\footnotesize $4\,\mathrm{m}$}}
    \put(197,115){{\footnotesize $3\,\mathrm{m}$}}
    \put(25,120){{\footnotesize member 1}}
    \put(138,120){{\footnotesize member 2}}
    \put(68,151){{\footnotesize $p_{x}$}}
    \put(104,174){{\footnotesize $p_{y}$}}
    }
  \end{picture}
  \medskip
  \caption{A two-bar truss. }
  \label{fig:m_2bar}
\end{figure}
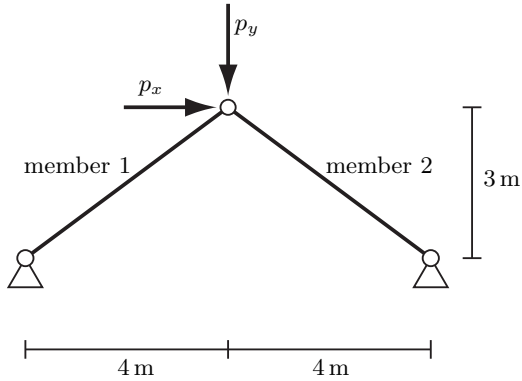

\begin{figure}[tp]
  \centering
  \includegraphics[scale=0.70]{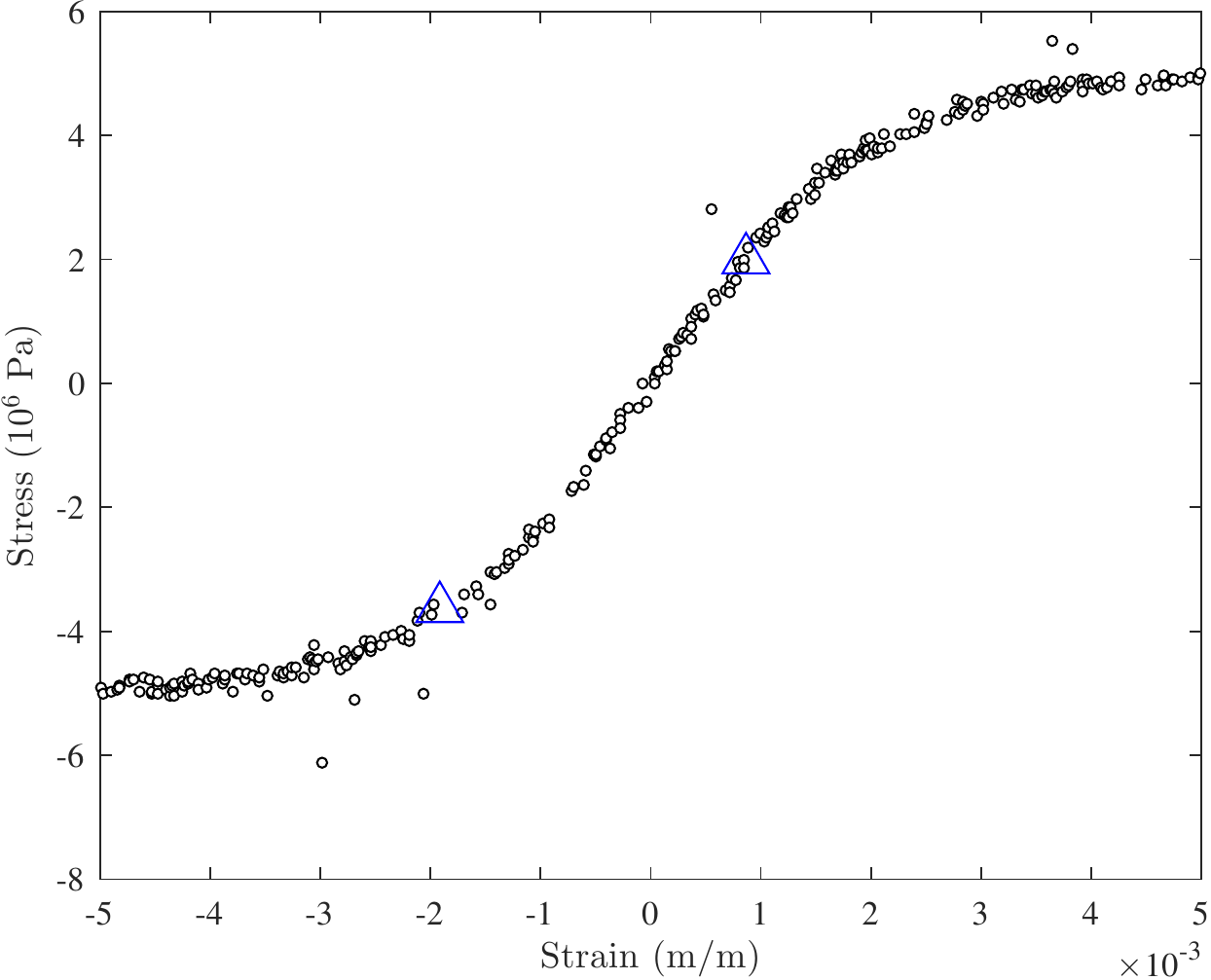}
  \caption{A material data set and the solution obtained by the proposed method. }
  \label{fig:ex_regress_constitutive}
\end{figure}

\begin{figure}[tp]
  \centering
  \includegraphics[scale=0.60]{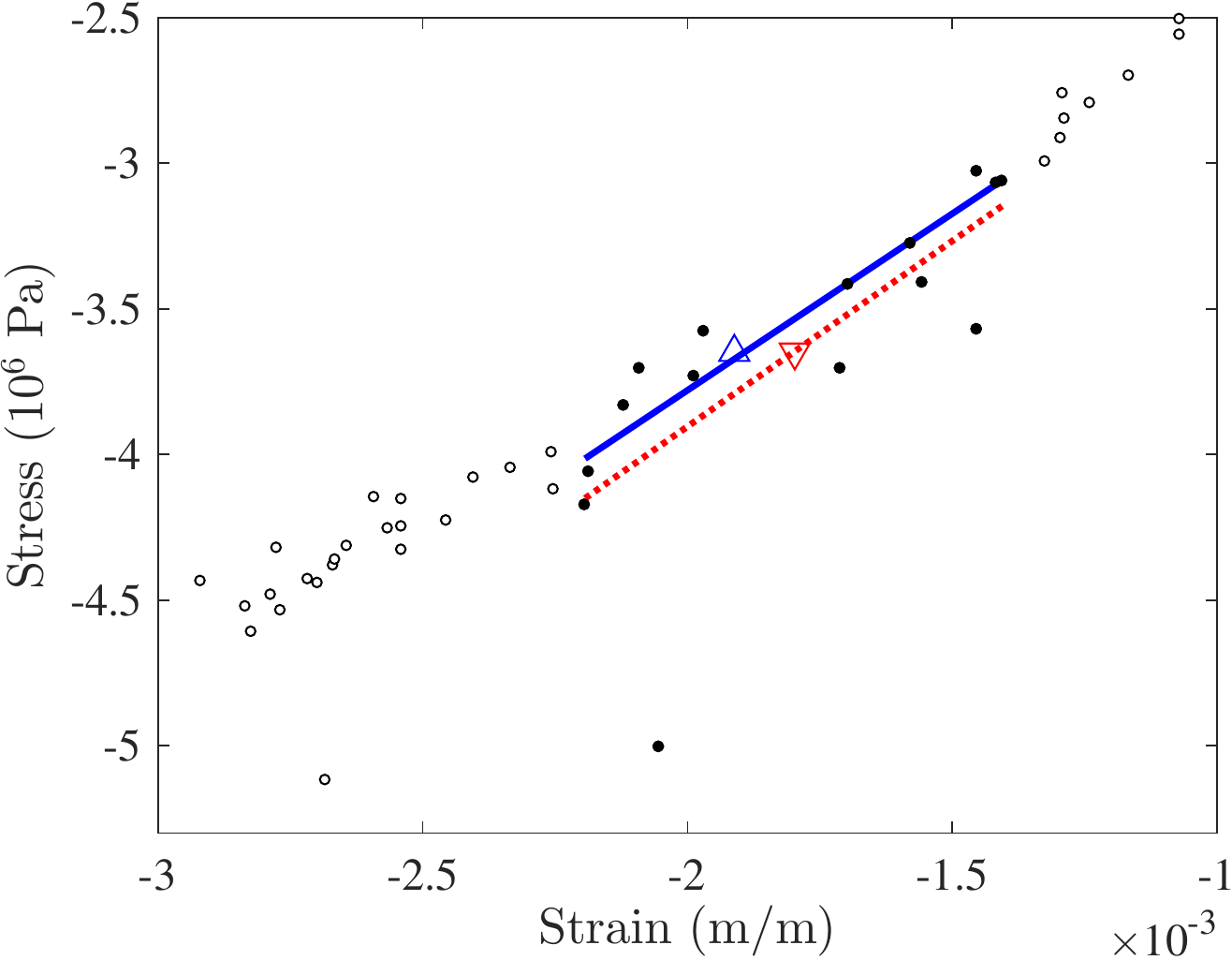}
  \caption{A closeup of \reffig{fig:ex_regress_constitutive}. 
  ``{\em solid line\/}'' The result obtained by the robust regression; 
  ``$\bullet$'' the data used for the robust regression; and 
  ``{\em dashed line\/}'' the result obtained by the least-squares regression. 
  }
  \label{fig:ex_regress_closeup}
\end{figure}

Consider the truss structure shown in \reffig{fig:m_2bar}. 
The two members consist of the same material. 
Suppose that the stress--strain relation of the material is 
characterized by the experimental data shown in 
\reffig{fig:ex_regress_constitutive}. 
Throughout the paper, we assume the elasticity. 
The cross-sectional area of each member is $2000\,\mathrm{mm^{2}}$. 
As for the external load, the horizontal force of 
$p_{x}=10\,\mathrm{kN}$ and the vertical force of $p_{y}=2\,\mathrm{kN}$ 
are applied at the top node. 

As the solution obtained by the method developed in this paper, 
the two triangles in \reffig{fig:ex_regress_constitutive} indicate 
the pairs of the stress and the strain of the two members. 
Here, the right triangle corresponds to member 1, and the left one 
corresponds to member 2. 
These stresses exactly satisfy the force-balance equation with the 
specified external load. 
Also, there exists a nodal displacement vector that exactly satisfies 
the compatibility relation with these strains.\footnote{
Since this example is a statically determinate truss, for any member 
strains there exists a compatible nodal displacement vector. 
What is meant to be explained here is that in the proposed method a 
solution is alway defined so as to satisfy the compatibility relation exactly. 
} 
Moreover, each pair of stress and strain seems to agree well with the 
data in \reffig{fig:ex_regress_constitutive}. 
\reffig{fig:ex_regress_closeup} shows a closeup of 
\reffig{fig:ex_regress_constitutive}. 

To obtain the solution above, we may consider that, in a certain sense, 
the data points that are far from the two triangles in 
\reffig{fig:ex_regress_constitutive} are unnecessary. 
This observation motivates us to develop a method that utilizes only 
the data points that are close to the final solution. 
In contrast, the conventional model-based methods in computational 
mechanics first assume an empirical single model describing the global 
behavior of the material. 
Then the parameters of the model is calibrated so that the model fits 
the material data. 
In the example in \reffig{fig:ex_regress_constitutive}, for instance, 
this model is required to fit the data points that are far from the 
solutions, as well as the near data points. 
However, in fact, fitting far data points is not necessary for agreement 
of the final solution with the data set. 
In contrast, a set of simple local models, which fit well only to some 
data points that are close to the solution, may probably yield a solution 
that is likely to reflect more the material data. 
This is because, from the bias-variance trade-off 
\cite{HTF09,Fla12}, a linear model has low variance and can have law 
bias for local fitting, while a global nonlinear model has high variance 
when we reduce its bias. 
This is an essential idea of the method presented in this paper. 

In \reffig{fig:ex_regress_closeup}, the solution indicated by 
``$\bigtriangleup$'' is computed by using only $15$ data points 
indicated by ``$\bullet$''. 
Namely, in the proposed method, we calibrate a linear model depicted by 
the solid line so as to fit well only these data, and use the calibrated 
linear model to perform the equilibrium analysis. 
Constructing such a local model from the data set is known as the 
$k$-nearest neighbor (kNN) local regression; the example in 
\reffig{fig:ex_regress_closeup} corresponds to $k=15$. 
See, e.g., \cite{Alt92,Cle79,CD88} 
for fundamentals of kNN local regression, 
and \cite{LOY17,PRFZ05,GGC08} for recent applications. 

Using only a few data points has, however, a disadvantage that the 
computational results can possibly be sensitive to outliers (or large 
errors). 
Indeed, if we use a conventional least-squares regression to calibrate 
a linear model, then we obtain the model depicted by the dotted line in 
\reffig{fig:ex_regress_closeup}. 
This model gives a lower stress level than the main locus of the data 
points, which shows the effect of one (or, possibly two) outlier(s) with 
small stress level(s). 
To avoid such undesired perturbation due to outliers, in this paper we 
adopt the robust regression. 
The solid line in \reffig{fig:ex_regress_closeup}, which is actually 
less affected by the outlier(s), is obtained by the robust regression. 

Because our local models are constructed by the linear regression, each model 
may be considered a parametric model. 
However, we do not know {\em a priori\/} which data points will be used 
to construct a local model. 
The set of data points to be used depends on the structural system, the 
external load, etc., and hence a set of the local models is adaptive. 
Therefore, in general, the set of the local models used in this method does 
not reflect the {\em global\/} behavior of the material. 
This is a significant difference of this method from the conventional 
model-based approach.

\section{Data-driven equilibrium analysis with robustness against noise and outliers}
\label{sec:method}

The method overviewed in section~\ref{sec:overview} is formally 
described in this section. 
In section~\ref{sec:method.regression}, we apply the local robust 
linear regression to  material data. 
In section~\ref{sec:method.equilibrium}, we present the overall 
framework of data-driven equilibrium analysis. 

\subsection{$k$-nearest neighbor local robust regression for material data}
\label{sec:method.regression}

As for the characterization of the material, suppose that the experimental 
data of the uniaxial strain, $\varepsilon$, and the uniaxial stress, 
$\sigma$, are given. 
The data set, denoted $D$, consists of pairs $(\varepsilon,\sigma)$ as 
\begin{align*}
  D = \{ (\check{\varepsilon}_{1}, \check{\sigma}_{1}), \dots, 
  (\check{\varepsilon}_{d}, \check{\sigma}_{d}) \} , 
\end{align*}
where $d$ is the number of data points. 
For simplicity, we assume that the structure consists of a single material. 
An example of a data set is shown in \reffig{fig:ex_regress_constitutive}. 
Define $E$ by 
\begin{align*}
  E = \{ \check{\varepsilon}_{1},\dots,\check{\varepsilon}_{d} \} , 
\end{align*}
which is the set of strain data points. 

For point $\varepsilon \in \Re$, 
let $N_{k}(\varepsilon) \subseteq E$ denote the set of $k$ data 
points that are closest to $\varepsilon$, i.e., 
\begin{align}
  N_{k}(\varepsilon) 
  = \argmin \Bigl\{ 
  \sum_{\check{\varepsilon} \in N}
  | \check{\varepsilon} - \varepsilon |
  \,\Bigm|
  N \subseteq E , \
  |N| = k 
  \Bigr\} . 
  \label{def:k.nearest.neighbor}
\end{align}
Namely, $N_{k}(\varepsilon)$ is the $k$-nearest neighbor of $\varepsilon$. 
We use $J_{k}(\varepsilon)$ to denote the set of indices of the 
data points in $N_{k}(\varepsilon)$, i.e., 
\begin{align*}
  J_{k}(\varepsilon) 
  = \{ j \in \{ 1,\dots,d \}
  \mid
  \check{\varepsilon}_{j} \in N_{k}(\varepsilon) \} . 
\end{align*}

Following the concept of the kNN local regression \cite{Cle79,CD88,Alt92}, 
we assume that, in the neighborhood of a given point 
$\tilde{\varepsilon} \in \Re$, the stress $\sigma$ is given by a 
deterministic linear function of $\varepsilon$ with additive noise as 
\begin{align}
  \sigma = w \varepsilon + v + \epsilon  , 
  \quad \forall \varepsilon \in \convex N_{k}(\tilde{\varepsilon}) , 
  \label{eq.regression}
\end{align}
where $w \in \Re$ and $v \in \Re$ are parameters, and 
$\varepsilon$ and $\sigma$ are considered the explanatory variable and 
the dependent variable, respectively. 
If $\epsilon$ is Gaussian noise with $0$ mean and 
$\mu^{2}$ variance (i.e., $\epsilon \sim \NC(0, \mu^{2})$), 
then the maximum likelihood estimation of the parameters of 
the regression model in \eqref{eq.regression} 
coincides with the least-squares 
approximation for the data points in $N_{k}(\tilde{\varepsilon})$ (and 
the corresponding stress data). 
This local least-squares method can be formulated as follows: 
\begin{align}
  \text{Minimize} 
  \quad
  \sum_{j \in J_{k}(\tilde{\varepsilon})} 
  (w \check{\varepsilon}_{j}  + v - \check{\sigma}_{j} )^{2} . 
  \label{eq.least.squares.1}
\end{align}
Here, $w$ and $v$ are variables to be optimized. 
Since we suppose that the data set $D$ includes some outliers, we replace 
the quadratic penalty function in \eqref{eq.least.squares.1} with a 
penalty function (or a loss function) for robust regression, 
denoted $\phi: \Re \to \Re$, as follows: 
\begin{align}
  \text{Minimize} 
  \quad
  \sum_{j \in J_{k}(\tilde{\varepsilon})} 
  \phi(w \check{\varepsilon}_{j}  + v - \check{\sigma}_{j} ) .  
  \label{eq.robust.regression.1}
\end{align}
An example of $\phi$ is the Huber penalty function \cite{MMY06,HTF09} 
\begin{align}
  \phi(t) = 
  \begin{dcases*}
    t^{2} 
    & if $|t| \le M$, \\
    M(2|t| - M) 
    & otherwise, 
  \end{dcases*}
  \label{def:Huber.function}
\end{align}
where $M > 0$ is a constant. 
The Huber penalty function works in 
the manner same as the least-squares penalty function for a residual 
smaller than $M$, and behaves as the $\ell_{1}$-norm penalty function 
for a residual larger than $M$. 
It is worth noting that problem \eqref{eq.robust.regression.1} is an 
unconstrained convex optimization problem \cite{BV04}. 

\begin{remark}
  In \eqref{def:k.nearest.neighbor}, to define the $k$-nearest neighbor, 
  we define the distance between two data points as the difference of 
  their strain values, $| \check{\varepsilon} - \varepsilon |$. 
  Other definition of the distance can certainly be adopted. 
  For instance, following the idea in \cite{KO16}, we may use 
  $\displaystyle \frac{1}{2} c_{e} (\check{\varepsilon} - \varepsilon)^{2} 
  + \frac{1}{2} \frac{1}{c_{e}} (\check{\sigma} - \sigma)^{2} $ 
  with positive constant $c_{e}$ instead of 
  $| \check{\varepsilon} - \varepsilon |$. 
  \finbox
\end{remark}

\subsection{Equilibrium analysis of truss structures}
\label{sec:method.equilibrium}

Consider a truss structure consisting of $m$ members. 
We use $n$ to denote the number of degrees of freedom of the nodal 
displacements. 

Let $\bi{u} \in \Re^{n}$ denote the vector of the nodal displacements. 
We use $\varepsilon_{i} \in \Re$ to denote the uniaxial strain of member $i$. 
Under the assumption of the small deformation, the compatibility 
relation can be written in the form 
\begin{align*}
  \varepsilon_{i}  = \bi{b}_{i}^{\top} \bi{u} , 
\end{align*}
where $\bi{b}_{i} \in \Re^{n}$ is a constant vector. 

For member $i$, let $a_{i}$ and $l_{i}$ denote the cross-sectional area 
and the undeformed member length, respectively, which are considered 
constants. 
We use $\sigma_{i} \in \Re$ to denote the uniaxial stress. 
The force-balance equation can be written as 
\begin{align*}
  \sum_{i=1}^{m} a_{i} l_{i} \sigma_{i} \bi{b}_{i}  = \bi{p} , 
\end{align*}
where $\bi{p} \in \Re^{n}$ is the vector of the nodal external forces. 

The relation between $\varepsilon$ and $\sigma$ can be retrieved from the 
experimental material data in the manner proposed in 
section~\ref{sec:method.regression}. 
Therefore, our data-driven solver predicts the nodal displacements at the 
equilibrium state as the solution to the following system: 
\begin{alignat}{2}
  \varepsilon_{i}  &= \bi{b}_{i}^{\top} \bi{u} , 
  &{\quad}& i=1,\dots,m , 
  \label{eq.equilibrium.truss.1} \\
  \sum_{i=1}^{m} a_{i} l_{i} \sigma_{i} \bi{b}_{i}  &= \bi{p} , 
  \label{eq.equilibrium.truss.2} \\
  \sigma_{i} 
  &= \hat{w}_{i} \varepsilon_{i} + \hat{v}_{i}  , 
  &{\quad}& i=1,\dots,m , 
  \label{eq.equilibrium.truss.3} \\
  (\hat{w}_{i}, \hat{v}_{i} ) 
  &= \argmin
  \Bigl\{
  \sum_{j \in J_{k}(\varepsilon_{i})} 
  \phi(w_{i} \check{\varepsilon}_{j}  + v_{i} - \check{\sigma}_{j} )
  \Bigm|
  (w_{i},v_{i}) \in \Re^{2}
  \Bigr\} , 
  &{\quad}& i=1,\dots,m . 
  \label{eq.equilibrium.truss.4}
\end{alignat}

It is worth noting that \eqref{eq.equilibrium.truss.1}, 
\eqref{eq.equilibrium.truss.2}, and \eqref{eq.equilibrium.truss.3} are 
linear equations. 
In contrast, it is not easy to deal with \eqref{eq.equilibrium.truss.4} 
directly, since its right-hand side involves unknown $\varepsilon_{i}$. 
This requires an iterative method. 
We adopt the following simple heuristic: 
Let $\varepsilon_{i}^{(l)}$ denote the incumbent solution of the strain of 
member $i$. 
For each $i=1,\dots,m$, we solve the optimization problem in 
\eqref{eq.equilibrium.truss.4} with $\varepsilon:=\varepsilon_{i}^{(l)}$ to 
obtain the optimal solution denoted by 
$(\hat{w}_{i}^{(l)}, \hat{v}_{i}^{(l)})$. 
With 
$(\hat{w}_{i}, \hat{v}_{i}):=(\hat{w}_{i}^{(l)}, \hat{v}_{i}^{(l)})$ in 
\eqref{eq.equilibrium.truss.3}, we solve the system of linear equations, 
\eqref{eq.equilibrium.truss.1}, \eqref{eq.equilibrium.truss.2}, and 
\eqref{eq.equilibrium.truss.3}, to find the solution 
$(\bi{u}^{(l+1)}, \bi{\varepsilon}^{(l+1)}, \bi{\sigma}^{{(l+1)}})$. 
We terminate the iteration if $\| \bi{u}^{(l+1)}-\bi{u}^{(l)} \|$ is small 
enough. 

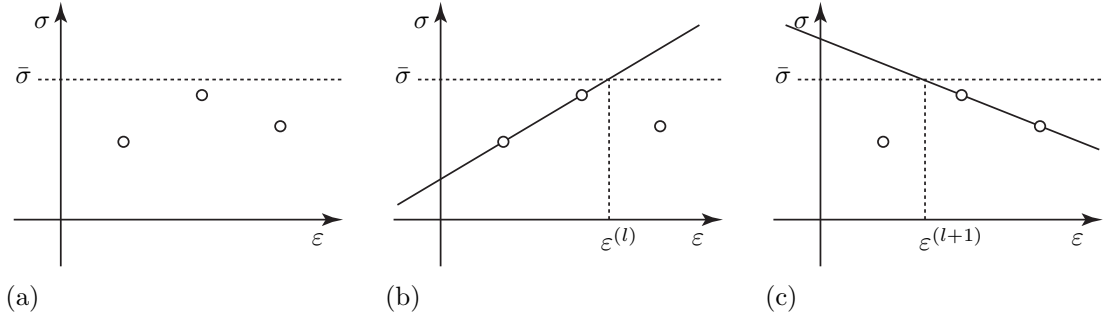
\begin{figure}[tp]
  \centering
  \begin{subfigure}[b]{0.32\textwidth}
    \centering
    \includegraphics[scale=0.65]{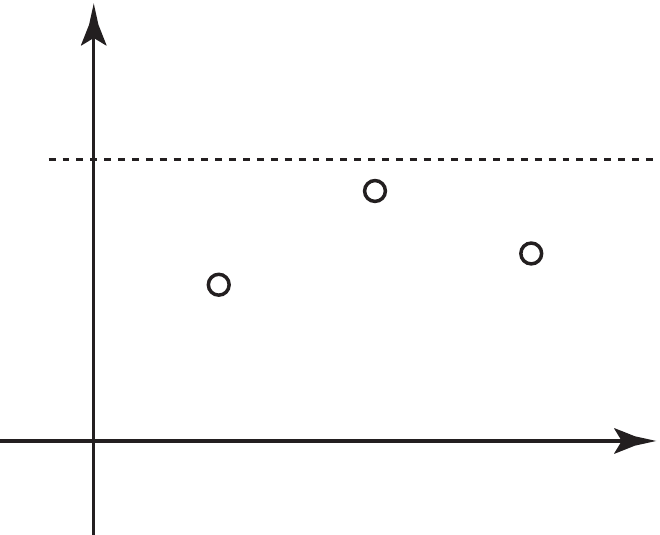}
    \begin{picture}(0,0)
      \put(-180,-10){
      \put(164,18){{\small $\varepsilon$}}
      \put(60,100){{\small $\sigma$}}
      \put(53,79){{\small $\bar{\sigma}$}}
      }
    \end{picture}
    \caption{}
    \label{fig:counter_ex1}
  \end{subfigure}
  \hfill
  \begin{subfigure}[b]{0.32\textwidth}
    \centering
    \includegraphics[scale=0.65]{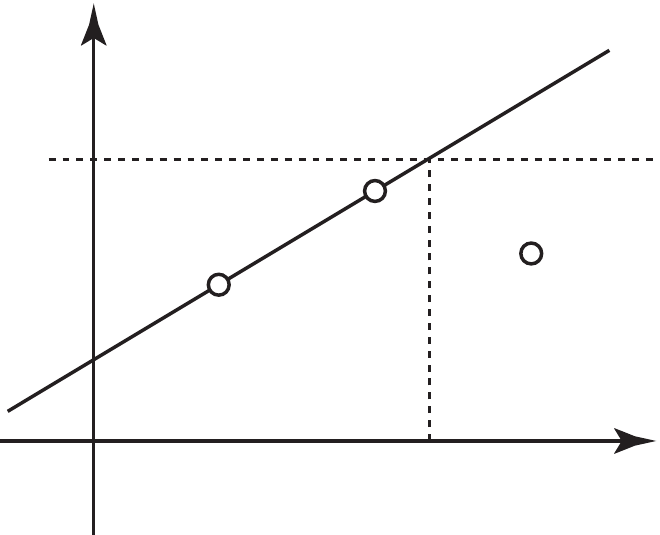}
    \begin{picture}(0,0)
      \put(-180,-10){
      \put(164,18){{\small $\varepsilon$}}
      \put(60,100){{\small $\sigma$}}
      \put(53,79){{\small $\bar{\sigma}$}}
      \put(130,16){{\small $\varepsilon^{(l)}$}}
      }
    \end{picture}
    \caption{}
    \label{fig:counter_ex2}
  \end{subfigure}
  \hfill
  \begin{subfigure}[b]{0.32\textwidth}
    \centering
    \includegraphics[scale=0.65]{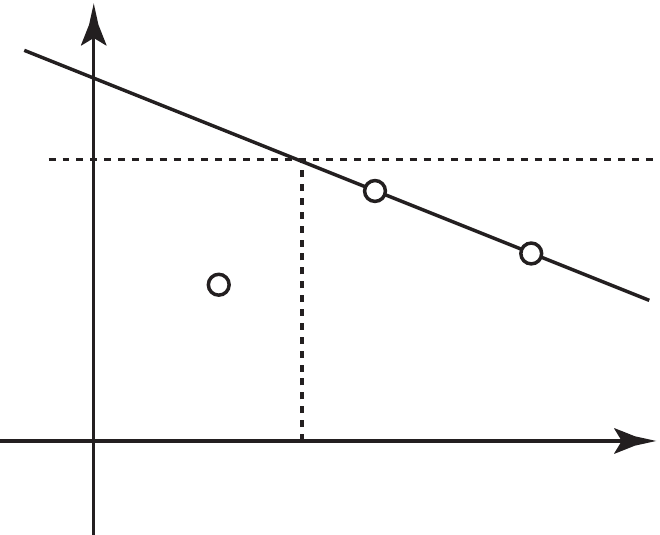}
    \begin{picture}(0,0)
      \put(-180,-10){
      \put(164,18){{\small $\varepsilon$}}
      \put(60,100){{\small $\sigma$}}
      \put(53,79){{\small $\bar{\sigma}$}}
      \put(107,16){{\small $\varepsilon^{(l+1)}$}}
      }
    \end{picture}
    \caption{}
    \label{fig:counter_ex3}
  \end{subfigure}
  \caption{An example such that the proposed problem formulation with 
  $k=2$ has no solution. 
  \subref{fig:counter_ex1} A material data set; 
  \subref{fig:counter_ex2} incumbent solution $\varepsilon^{(l)}$; and 
  \subref{fig:counter_ex3} solution $\varepsilon^{(l+1)}$ obtained from 
  the local linear regression for $N_{2}(\varepsilon^{(l)})$. 
  }
  \label{fig:counter_ex}
\end{figure}

\begin{remark}\label{rem:existence.of.solution}
  Since this paper is intended to be the first attempt to make use of the local 
  robust regression in the data-driven computational mechanics, the 
  proposed method admits of improvement. 
  Particularly, a drawback is that the system \eqref{eq.equilibrium.truss.1}, 
  \eqref{eq.equilibrium.truss.2}, \eqref{eq.equilibrium.truss.3}, and 
  \eqref{eq.equilibrium.truss.4} may possibly have no solution. 
  For instance, consider the material data set in \reffig{fig:counter_ex1}. 
  Suppose that a statically determinate truss is of interest, and the 
  stress of one of members, satisfying the force-balance equation, is $\bar{\sigma}$. 
  We set $k=2$. 
  Suppose that the 2-nearest neighbor of the incumbent solution consists 
  of the left two data points. 
  Then the local linear regression results in the solid line in 
  \reffig{fig:counter_ex2}, and the corresponding strain is obtained as 
  $\varepsilon^{(l)}$. 
  In the next iteration, $N_{2}(\varepsilon^{(l)})$ consists of the two 
  right data points. 
  Therefore, the local linear regression results in the solid line in 
  \reffig{fig:counter_ex3}. 
  Accordingly, the strain is updated to $\varepsilon^{(l+1)}$. 
  Since $N_{2}(\varepsilon^{(l+1)})$ consists of the two left data 
  points, we go back to the situation in \reffig{fig:counter_ex2}. 
  Thus, the system \eqref{eq.equilibrium.truss.1}, 
  \eqref{eq.equilibrium.truss.2}, \eqref{eq.equilibrium.truss.3}, and 
  \eqref{eq.equilibrium.truss.4} has no solution for the example in 
  \reffig{fig:counter_ex1} with $k=2$. 
  We leave resolution of this issue as future work, and attention of 
  this paper is focused on the effectiveness of the presented method in 
  avoiding the influence of 
  noise and outliers when a solution is successfully found. 
  In the numerical experiments in section~\ref{sec:ex}, a solution is 
  found for almost all problem instances. 
  In the case that the proposed method does not converge, we increase 
  the value of $k$ as a heuristic way of coping with this issue. 
  \finbox
\end{remark}

\section{Numerical experiments}
\label{sec:ex}

The method proposed in section~\ref{sec:method} was implemented 
on Matlab ver.~9.0.0. 
A Matlab bult-in function \texttt{robustfit} was used to solve the 
robust regression with the Huber penalty function in 
\eqref{eq.equilibrium.truss.4}. 
In section~\ref{sec:ex.1}, we compare the proposed method with the 
existing data-driven solver proposed in \cite{KO16}. 
In section~\ref{sec:ex.2}, we perform the comparison with the method 
using the least-squares regression, which is not robust against the 
presence of outliers.

\subsection{Example (I): One-bar truss}
\label{sec:ex.1}

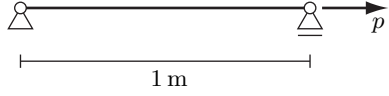
\begin{figure}[tp]
  \centering
  \includegraphics[scale=0.50]{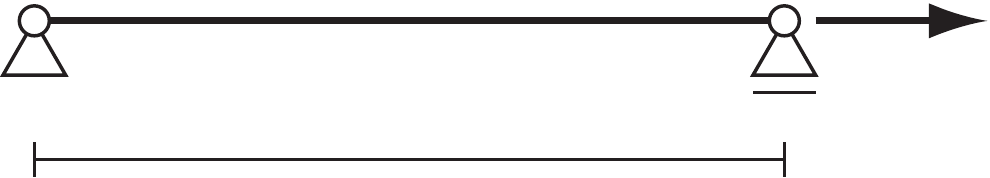}
  \begin{picture}(0,0)
    \put(-200,-50){
    \put(190,65){{\footnotesize $p$}}
    \put(107,43){{\footnotesize $1\,\mathrm{m}$}}
    }
  \end{picture}
  \medskip
  \caption{A one-bar truss. }
  \label{fig:m_1bar}
\end{figure}

\begin{figure}[tp]
  \centering
  \includegraphics[scale=0.60]{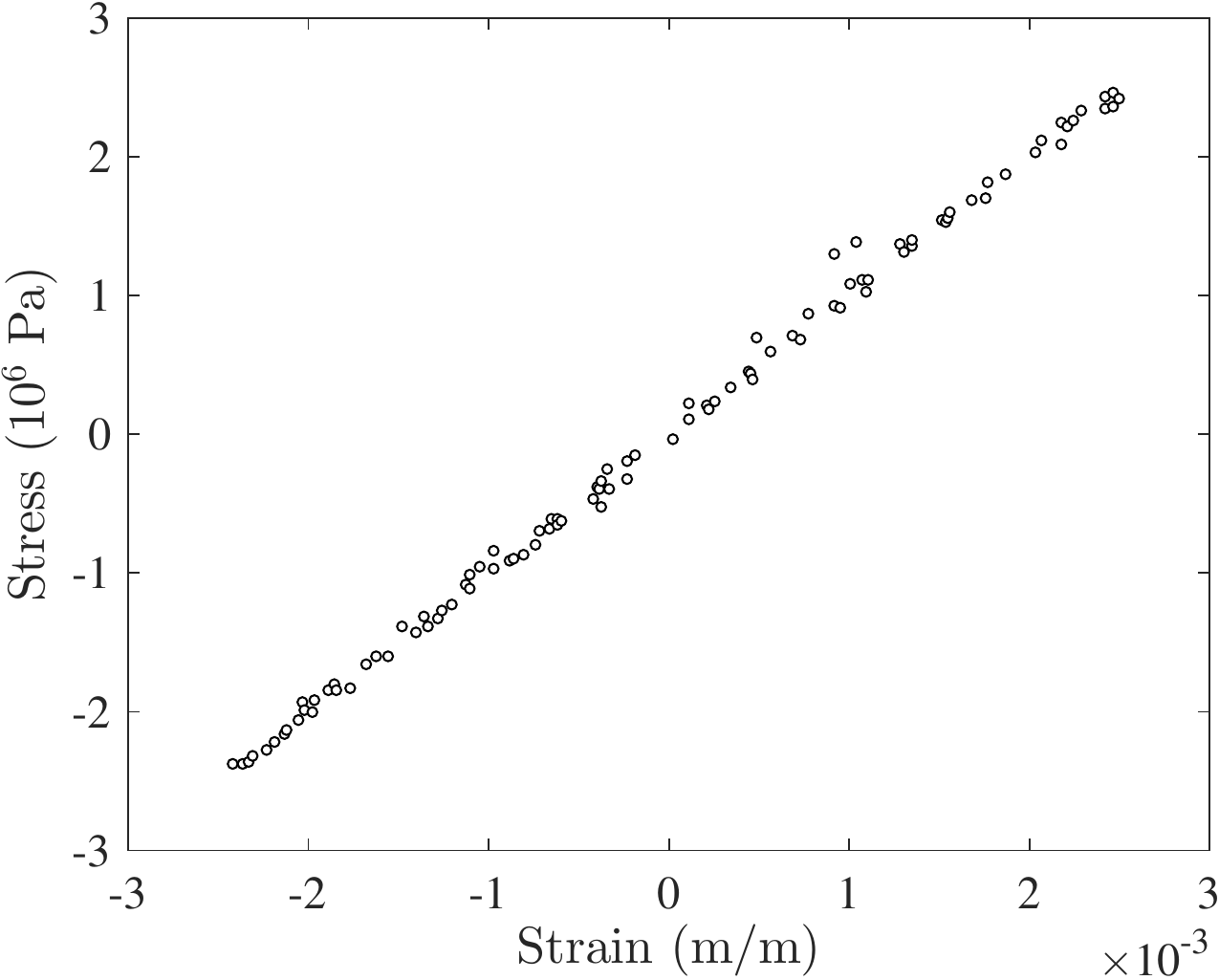}
  \caption{The material data set used for example (I). 
  }
  \label{fig:one_bar_constitutive}
\end{figure}

\begin{figure}[tp]
  \centering
  \begin{subfigure}[b]{0.63\textwidth}
    \centering
    \includegraphics[scale=0.60]{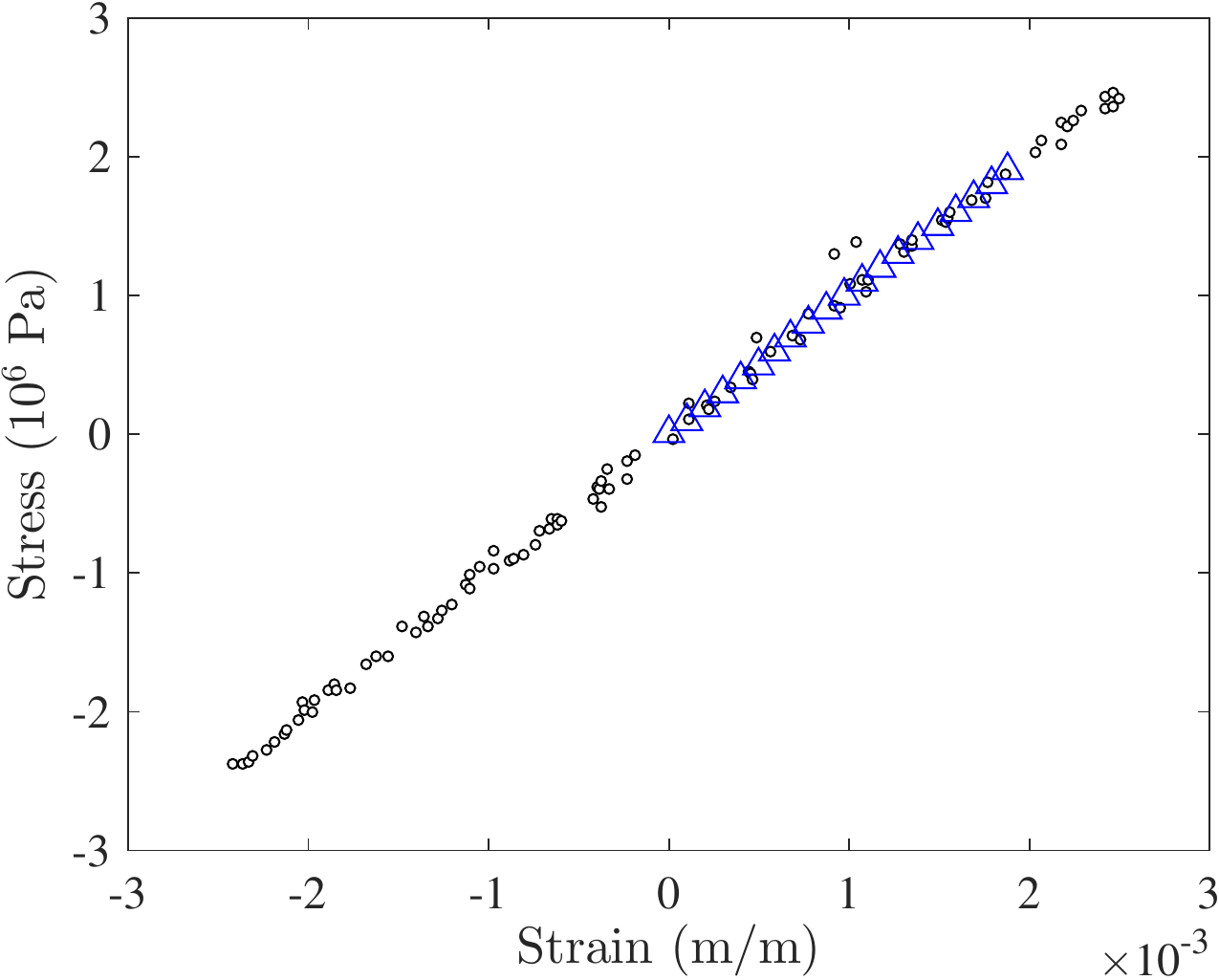}
    \caption{}
    \label{fig:one_bar_robust_constitutive}
  \end{subfigure}
  \hfill
  \begin{subfigure}[b]{0.36\textwidth}
    \centering
    \includegraphics[scale=0.40]{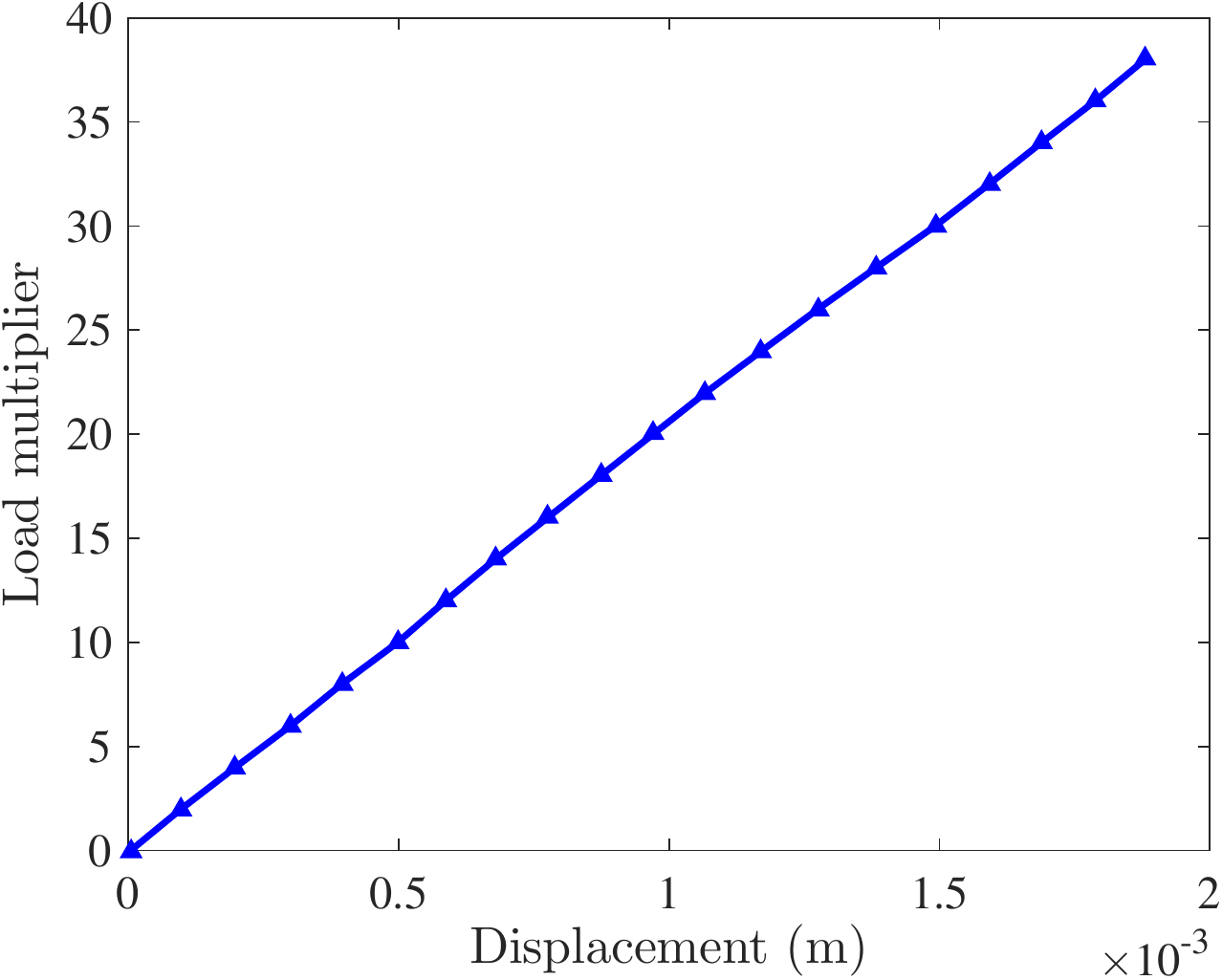}
    \caption{}
    \label{fig:one_bar_robust_load}
  \end{subfigure}
  \caption{The results obtained by the proposed method (example (I)). 
  \subref{fig:one_bar_robust_constitutive} The plot of the 
  solutions on the space of $(\varepsilon,\sigma)$; and 
  \subref{fig:one_bar_robust_load} the equilibrium path. 
  }
  \label{fig:one_bar_robust}
\end{figure}

\begin{figure}[tp]
  \centering
  \begin{subfigure}[b]{0.63\textwidth}
    \centering
    \includegraphics[scale=0.60]{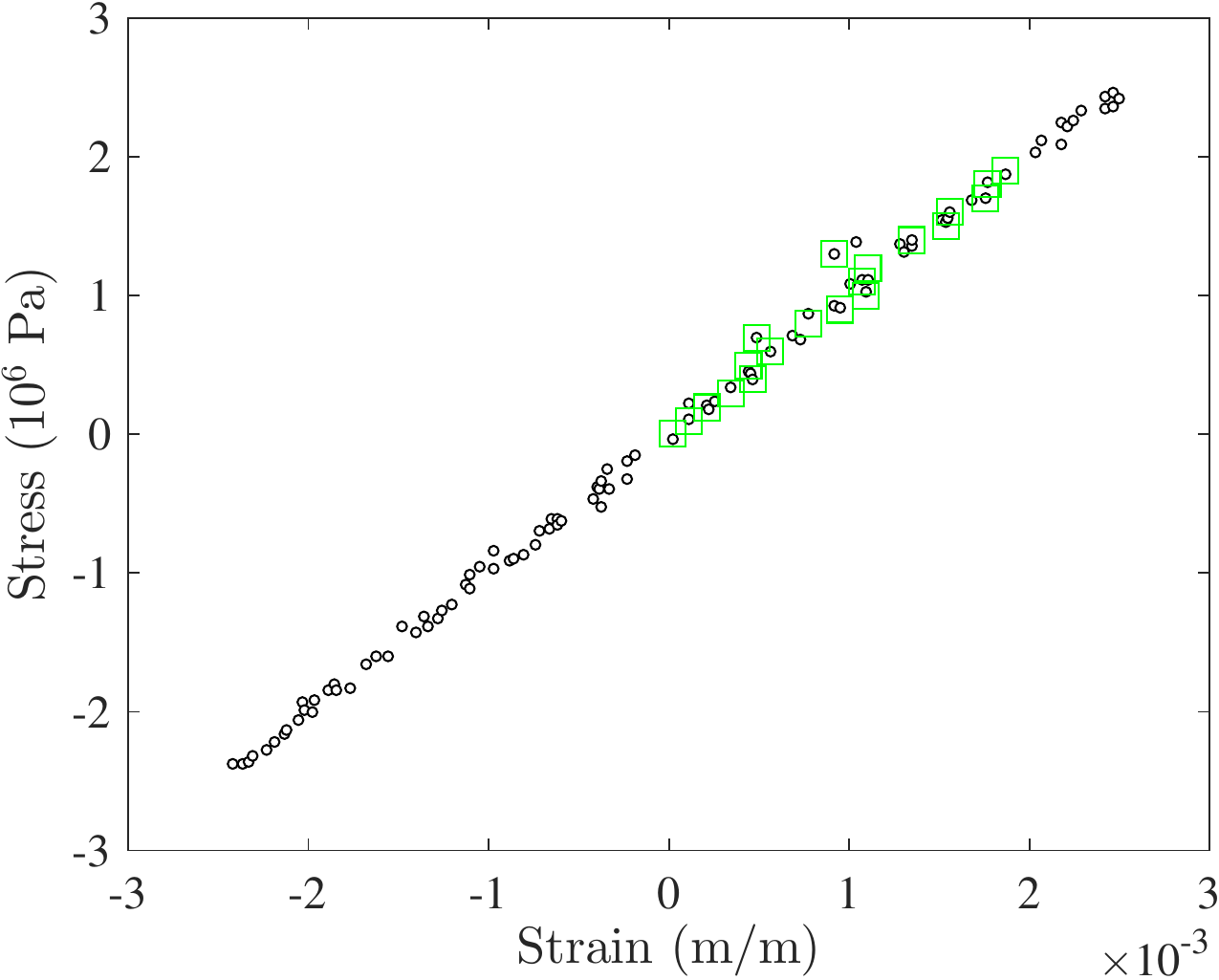}
    \caption{}
    \label{fig:one_bar_min_constitutive}
  \end{subfigure}
  \hfill
  \begin{subfigure}[b]{0.36\textwidth}
    \centering
    \includegraphics[scale=0.40]{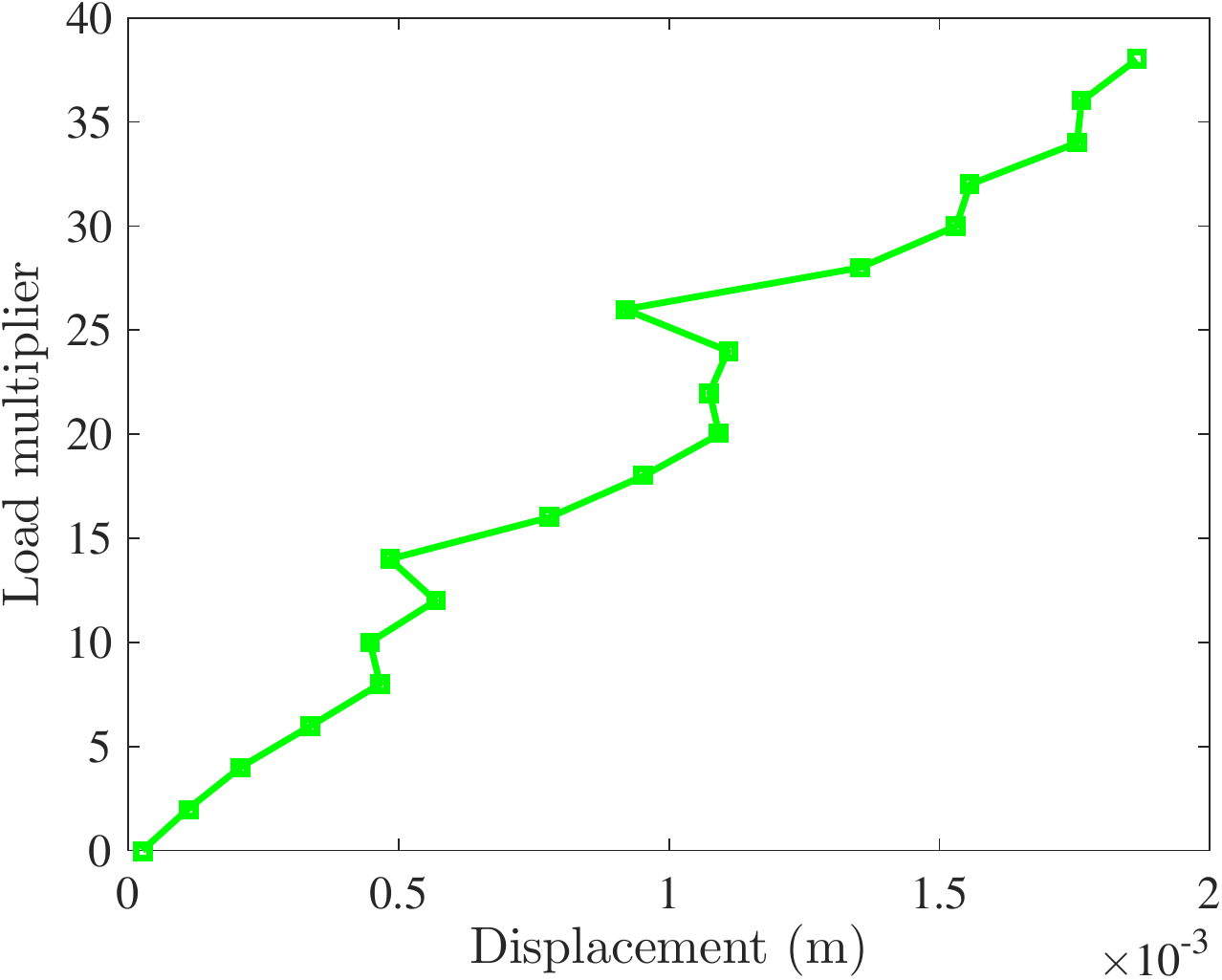}
    \caption{}
    \label{fig:one_bar_min_load}
  \end{subfigure}
  \caption{The results obtained by the method in \cite{KO16} (example (I)). 
  \subref{fig:one_bar_min_constitutive} The plot of the solutions on the 
  space of $(\varepsilon,\sigma)$; and 
  \subref{fig:one_bar_min_load} the equilibrium path. 
  }
  \label{fig:one_bar_min}
\end{figure}

As for the simplest example, consider a bar in \reffig{fig:m_1bar}. 
The bar cross-sectional area is $a=200\,\mathrm{mm^{2}}$. 
The external load is given as $p = \lambda \bar{p}$ with 
$\bar{p}=10\,\mathrm{N}$, where $\lambda$ is the load multiplier. 

In the proposed method, the size of the neighborhood is set to $k=15$. 
The \texttt{tune} parameter of the Matlab bult-in function 
\texttt{robustfit} is set to $10^{-3}$, which was determined by 
preliminary numerical experiments. 

\reffig{fig:one_bar_constitutive} shows the material data set, which 
consists of $100$ data points. 
This data set might be interpreted as the linear elasticity with 
additive noise. 

The computational results of the proposed method are shown in 
\reffig{fig:one_bar_robust}. 
The equilibrium states were computed for 
$\lambda = 0$, $2$, $4,\dots,38$. 
The triangles in \reffig{fig:one_bar_robust_constitutive} indicate the 
values of $(\varepsilon,\sigma)$ at these equilibrium states. 
\reffig{fig:one_bar_robust_load} plots the variation of the displacement 
with respect to the load multiplier. 
The linearity of the structural behavior is drawn out despite the 
presence of noise in the material data. 

\reffig{fig:one_bar_min} shows the results obtained by the method 
proposed in \cite{KO16}. 
In contrast to \reffig{fig:one_bar_robust_load}, it is clearly seen that 
the equilibrium path in \reffig{fig:one_bar_min_load} is strongly 
affected by noise. 
It is worth noting that, both in the methods proposed in \cite{KO16} and 
in this paper, the member stresses of a statically determinate truss are 
uniquely determined from the force-balance equation (in this example, 
the stress is determined as $\lambda\bar{p}/a$). 
Moreover, in this example the compatibility relation is automatically 
satisfied for any strain. 
Therefore, the method in \cite{KO16} adopt the data point the stress 
value of which is closest to $\lambda\bar{p}/a$. 
The displacement is then determined from the strain value of that data 
point. 
Use of information of only a single data point makes the method in 
\cite{KO16} sensitive to noise.

\subsection{Example (II): 27-bar truss}
\label{sec:ex.2}

\begin{figure}[tp]
  \centering
  \includegraphics[scale=0.60]{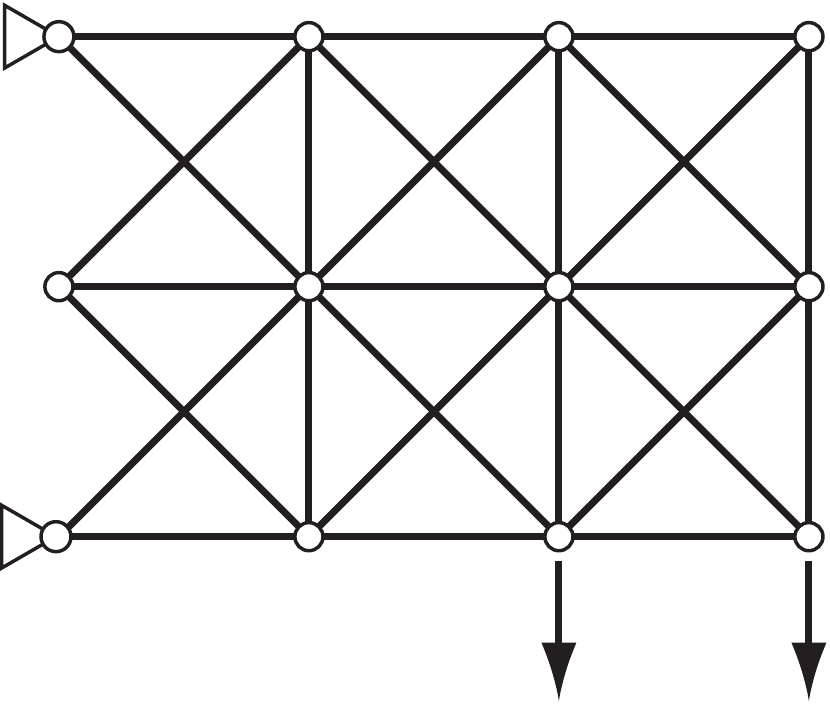}
  \caption{A 27-bar truss. }
  \label{fig:m_27bar}
\end{figure}

\begin{figure}[tp]
  \centering
  \begin{subfigure}[b]{0.47\textwidth}
    \centering
    \includegraphics[scale=0.50]{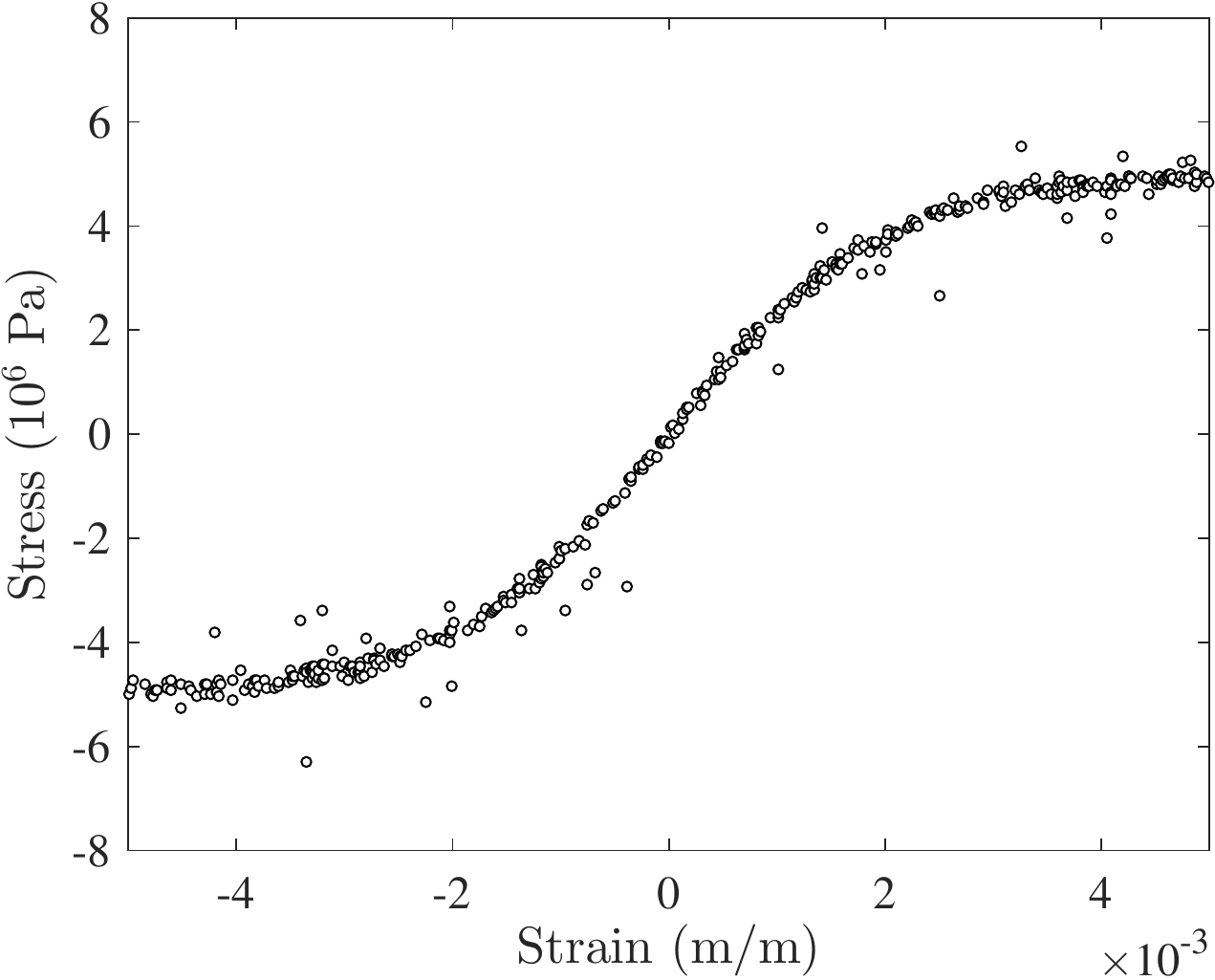}
    \caption{}
    \label{fig:27bar_data_1}
  \end{subfigure}
  \hfill
  \begin{subfigure}[b]{0.47\textwidth}
    \centering
    \includegraphics[scale=0.50]{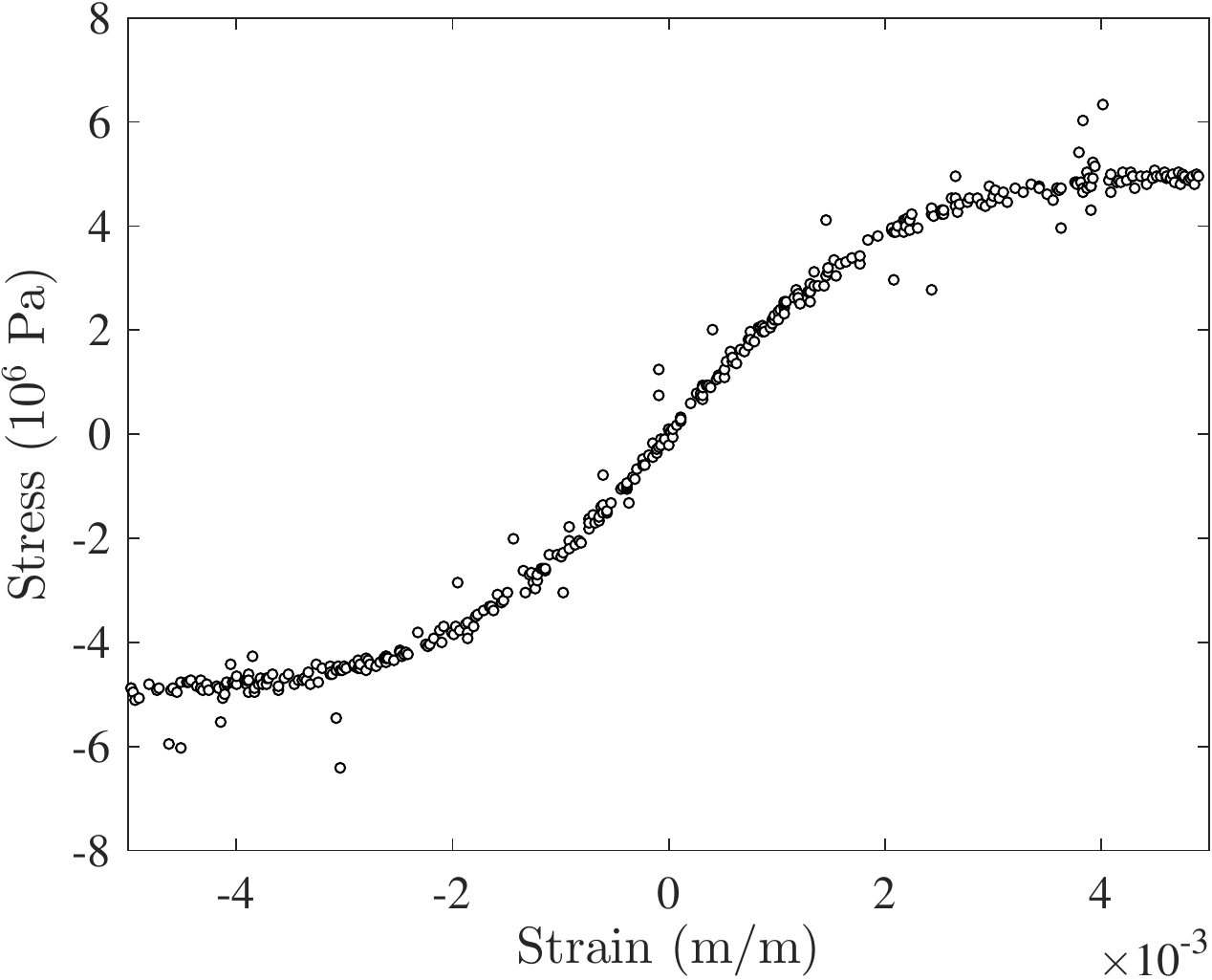}
    \caption{}
    \label{fig:27bar_data_3}
  \end{subfigure}
  \par\medskip
  \begin{subfigure}[b]{0.47\textwidth}
    \centering
    \includegraphics[scale=0.50]{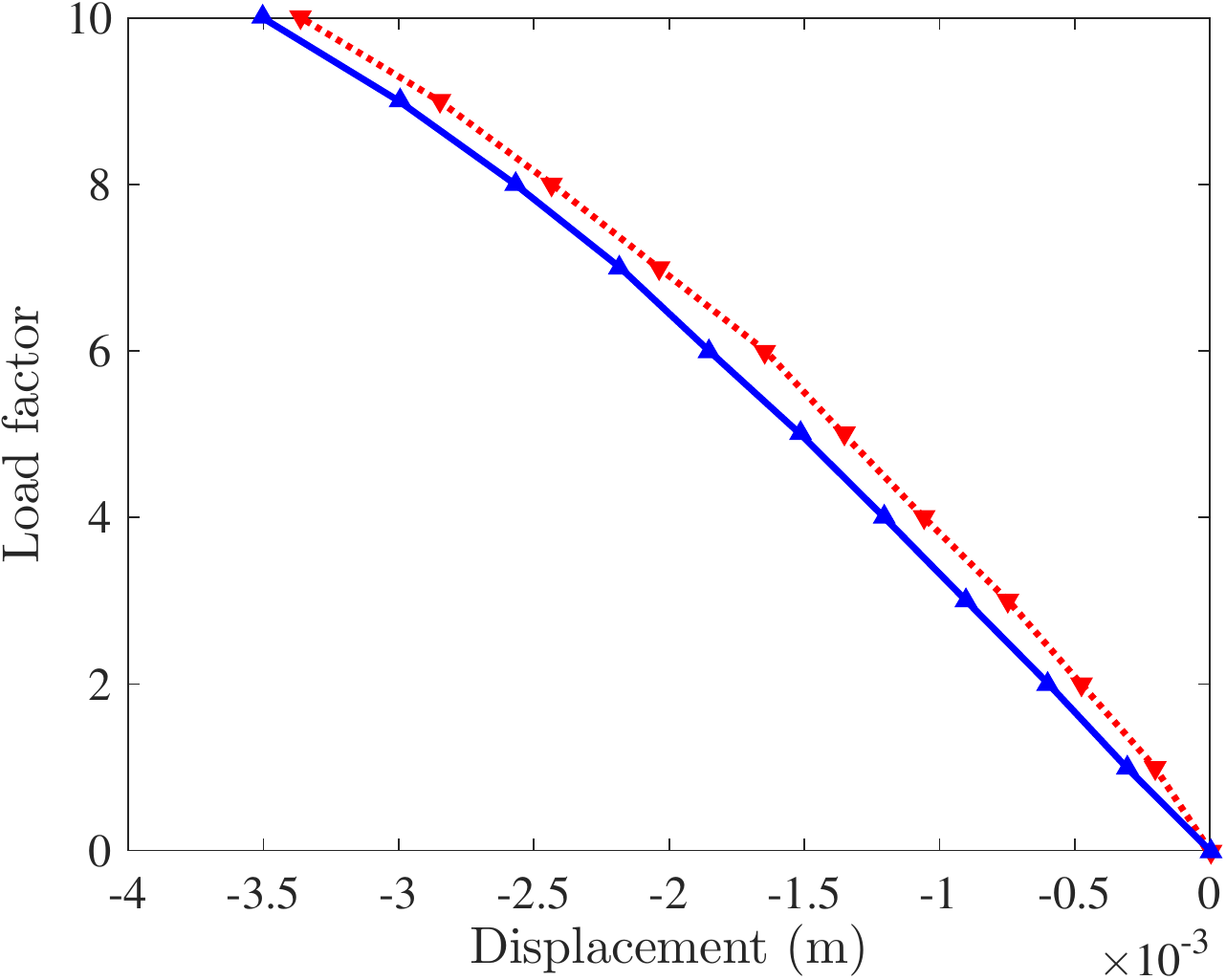}
    \caption{}
    \label{fig:27bar_path_1}
  \end{subfigure}
  \hfill
  \begin{subfigure}[b]{0.47\textwidth}
    \centering
    \includegraphics[scale=0.50]{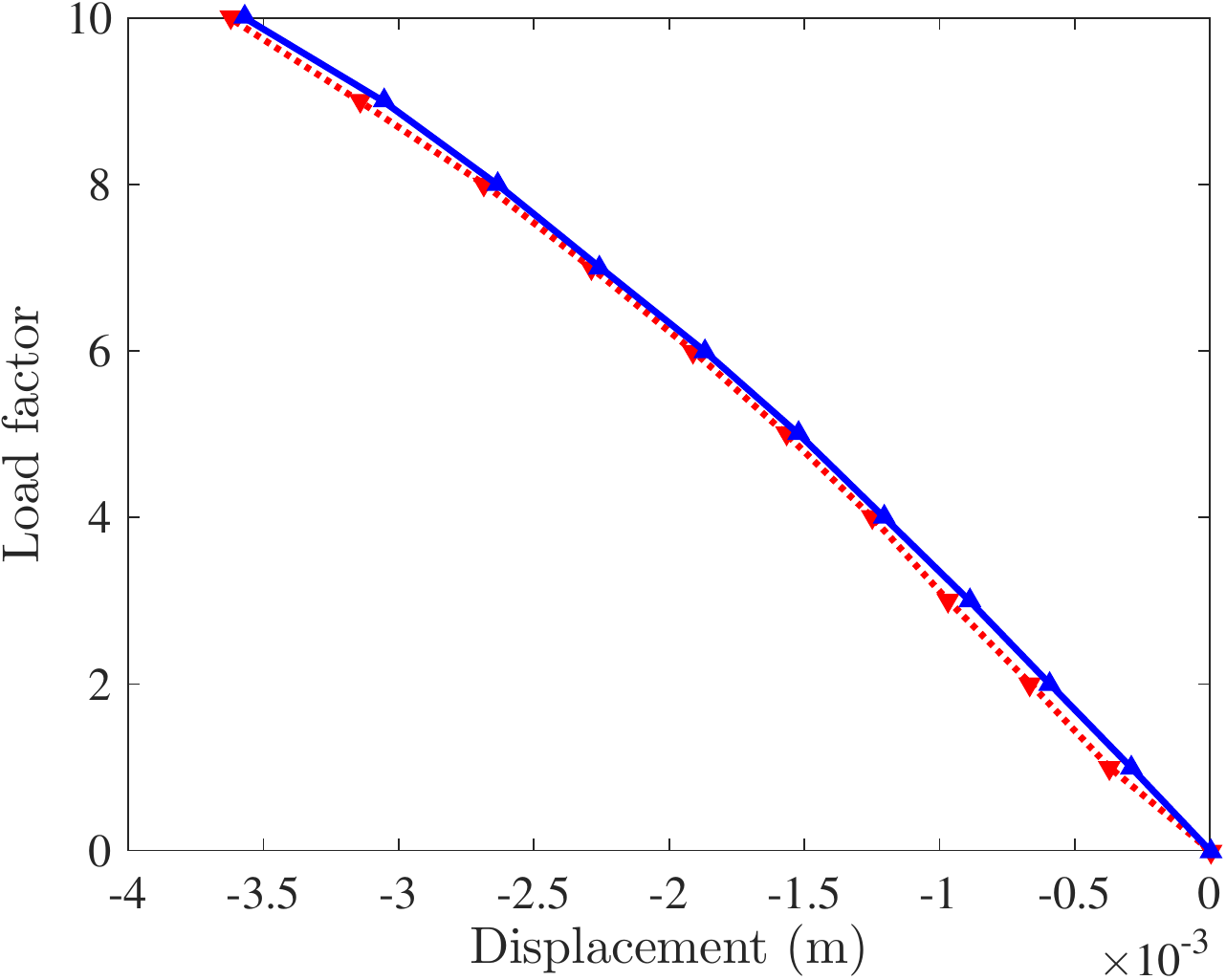}
    \caption{}
    \label{fig:27bar_path_3}
  \end{subfigure}
  \caption{Two typical computational results of example (II). 
  \subref{fig:27bar_data_1}, \subref{fig:27bar_data_3} 
  The material data sets; and 
  \subref{fig:27bar_path_1}, \subref{fig:27bar_path_3} 
  the obtained equilibrium paths. 
  ``{\em solid line\/}'' The result of the proposed method; and 
  ``{\em dotted line\/}'' the result of the least-squares regression. 
  }
  \label{fig:27bar_equilibrium}
\end{figure}

\begin{figure}[tp]
  \centering
  \begin{subfigure}[b]{0.47\textwidth}
    \centering
    \includegraphics[scale=0.50]{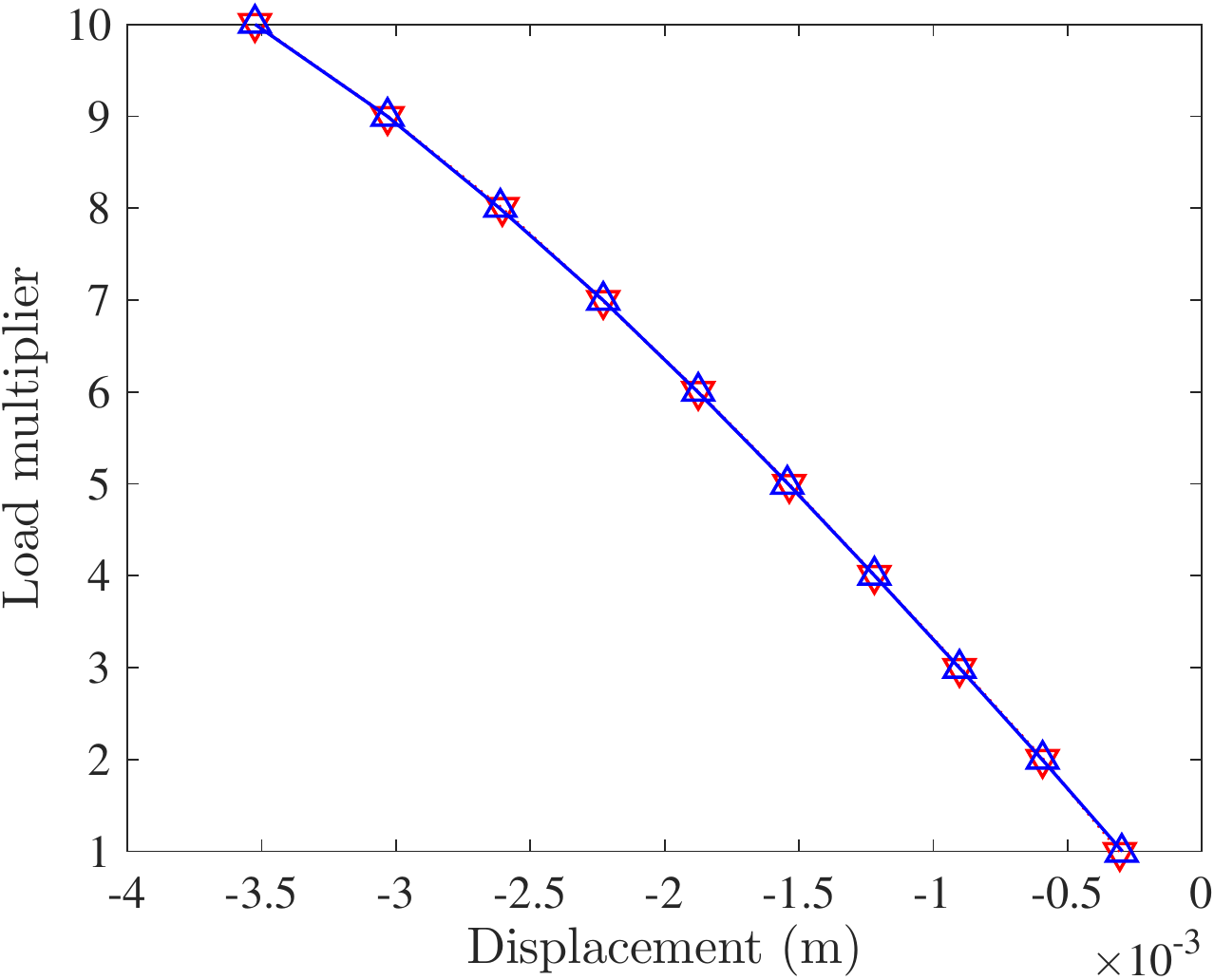}
    \caption{}
    \label{fig:27bar_mean_path}
  \end{subfigure}
  \hfill
  \begin{subfigure}[b]{0.47\textwidth}
    \centering
    \includegraphics[scale=0.50]{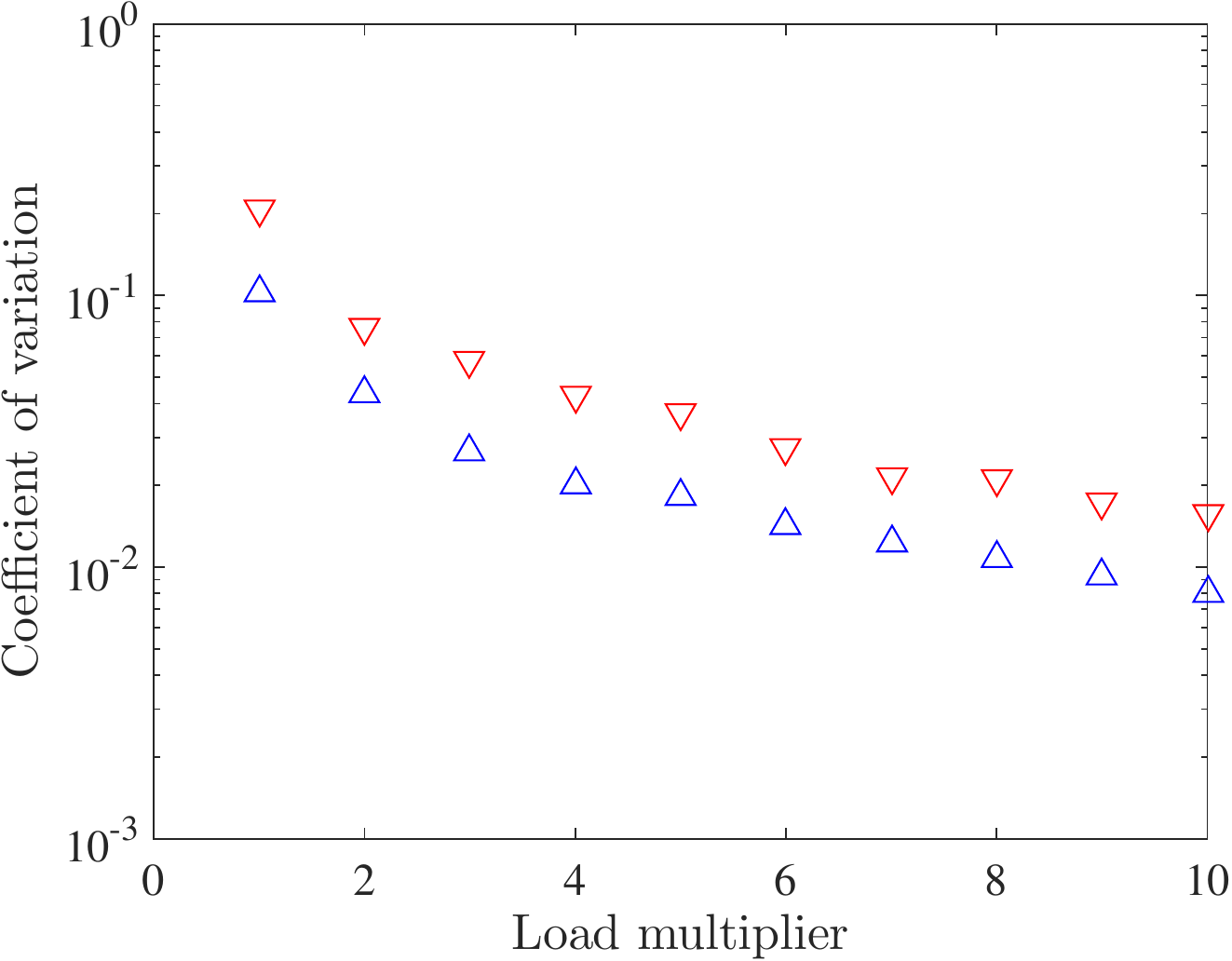}
    \caption{}
    \label{fig:27bar_variation_path}
  \end{subfigure}
  \caption{The statistics of the computational results of example (II). 
  \subref{fig:27bar_mean_path} The means; and 
  \subref{fig:27bar_variation_path} the coefficients of variation. 
  ``$\bigtriangleup$'' The proposed method; and 
  ``$\bigtriangledown$'' the least-squares regression. 
  }
  \label{fig:27bar_statistics}
\end{figure}

Consider the planar truss shown in \reffig{fig:m_27bar}. 
The undeformed length of each of the vertical and horizontal 
members is $1\,\mathrm{m}$. 
The cross-sectional area of every member is $200\,\mathrm{mm^{2}}$. 
As for the external load, $\bi{p}=\lambda\bar{\bi{p}}$, vertical 
downward force of $50\lambda$ in $\mathrm{N}$ are applied at the bottom 
two nodes as shown in \reffig{fig:m_27bar}. 

In the proposed method, the size of neighborhood is set to $k=40$. 
The \texttt{tune} parameter of \texttt{robustfit} is set to $10^{-4}$. 
For comparison, the method using the least-squares regression in 
\eqref{eq.least.squares.1} is also examined. 
A Matlab bult-in function \texttt{regress} was used for this purpose. 

\reffig{fig:27bar_data_1} and \reffig{fig:27bar_data_3} show 
typical data sets, each of which consists of $400$ data points. 
The computational results for these data sets are shown in 
\reffig{fig:27bar_path_1} and \reffig{fig:27bar_path_3}, respectively. 
The solid line depicts the equilibrium path (the variation of the 
horizontal displacement of the bottom rightmost node with respect to the 
load multiplier) obtained by the proposed method, while the dotted line 
depicts the equilibrium path obtained by using the least-squares 
regression. 
It is observed in \reffig{fig:27bar_path_1} that 
that the magnitudes of the displacements obtained by the least-squares 
regression are smaller than the ones obtained by the proposed method. 
In contrast, in \reffig{fig:27bar_path_3}, the magnitudes of the displacements obtained by 
the least-squares regression are larger than the ones obtained by the 
proposed method. 
Thus, the results of the least-squares regression are more strongly 
affected by outliers, and seem to have larger variation depending on 
given data sets. 
This observation is verified in \reffig{fig:27bar_statistics}, as 
explained in the following. 

We generated a data set $D$ consisting of $400$ observations as follows. 
The data set $E$ of the strain is drawn from the uniform distribution on 
interval $[-5\times 10^{-3}, 5\times 10^{-3}]$. 
From among them, we randomly choose $360$ data points, and set the 
stress values as 
\begin{align*}
  \sigma = \tilde{\sigma} + 0.1\epsilon_{1}
\end{align*}
with 
\begin{align*}
  \tilde{\sigma} &= \frac{10^{6}}{1 + \exp(-10^{3} \cdot \varepsilon)} 
\end{align*}
in $\mathrm{Pa}$ and the noise $\epsilon_{1} \sim \NC(0,1)$. 
The remaining $40$ data points are possibly outliers, the stress values 
of which are set as 
\begin{align*}
  \sigma = \tilde{\sigma} + 0.1\epsilon_{1} + 0.8\epsilon_{2} 
\end{align*}
with $\epsilon_{2} \sim \NC(0,1)$. 
One data set is generated in this way. 
We generated $100$ data sets independently, for each of which we 
performed the equilibrium analysis with the proposed method, as well as 
the method with the least-squares regression. 
The computational results are shown in \reffig{fig:27bar_statistics}. 
\reffig{fig:27bar_mean_path} plots the means of the obtained 
displacements. 
It is observed that the results of the two methods agree well with each 
other. 
\reffig{fig:27bar_variation_path} shows (the absolute values of) 
coefficients of variation. 
It is observed that the coefficients of variation are reduced by using 
the robust regression instead of the least-squares regression. 
This illustrates that the proposed method is robust against the presence 
of outliers. 
It is worth noting that this robustness stems from the definition of 
Huber penalty function in \eqref{def:Huber.function}, where we can see 
that, for large $t$, the increase of $\phi(t)$ with respect to the 
increase of $t$ is less drastic than the increase of the quadratic 
penalty function used in the least-squares regression. 

It is worth noting that additional computational cost required by the 
proposed method, compared with a conventional method for problems with 
material nonlinearity, is not very large. 
In fact, a standard method for the equilibrium analysis with material 
nonlinearity may be application of the Newton--Raphson method, which 
usually takes several iterations before convergence. 
In the numerical experiments above, the proposed method with the robust 
regression requires $3.6$ iterations in average to obtain the results 
shown in \reffig{fig:27bar_path_3}, and the one with the least-squares 
regression requires $4.9$ iterations in average. 
Additional computational task of the proposed method is sorting the 
material data points for identifying $J_{k}(\varepsilon_{i}^{(l)})$ and 
performing the robust regression for finding 
$(\hat{w}_{i}^{(l)},\hat{v}_{i}^{(l)})$. 
Among then, the sorting can be carried out efficiently with time 
complexity $O(d\log d)$. 
Also, the robust regression problem considered in this paper is recast 
as a convex quadratic problem, which can be solved within 
the polynomial time of $k$ \cite{BV04}, where $k \ll d$. 
Thus, the two additional computational tasks of the proposed method can 
be performed efficiently. 

\section{Conclusions}
\label{sec:conclude}

In this paper, we have presented a simple heuristic for data-driven 
computational elasticity with data involving noise and outliers. 
In contrast to the conventional methods that assume a single model 
describing the entire stress--strain relationship, the proposed method 
selects some data points adaptively and construct a set of local models 
each of which corresponds to the material property of one structural element. 
The material data points that are not close to the final solution are 
ignored in the proposed method, and, hence, the set of local models can 
reflect the data more directly compared with a single entire model. 

The existing data-driven approach to computational elasticity 
\cite{KO16} uses information of a single data point for each 
structural element. 
To avoid the influence of noise included the data, the proposed method 
uses information of the $k$-nearest neighbor to construct a local model. 
Since only quite a small number of data points are used for 
constructing each local model, we adopt the robust regression to make 
the proposed method robust against the presence of outliers. 
Furthermore, the method is simple, easy to understand, 
and easy to implement. 

This paper has been intended to be the first attempt to develop a robust 
data-driven approach to computational mechanics. 
Much remains to be explored. 
For example, as already mentioned in \refrem{rem:existence.of.solution}, 
a solution does not necessarily exist to the proposed formulation, 
although in the numerical experiments a solution was found successfully 
for almost all instances. 
Also, it is not clear how one can choose an appropriate size of the 
neighborhood, $k$. 
An adaptive scheme to adjust the value of $k$ might be able to be 
explored. 
Extensions of the proposed method to inelastic problems might be challenging. 
In contrast, it seems to be straightforward 
to apply the proposed method to large deformation problems. 
Although we have restricted to truss structures for simple 
presentation, the proposed method can be applied to frame structures in 
a straightforward manner. 
Extension to continua remains to be studied.

\paragraph{Acknowledgments}

This work is partially supported by 
JSPS KAKENHI 17K06633.

\end{document}